# NO ARBITRAGE WITHOUT SEMIMARTINGALES

By Robert A. Jarrow, Philip Protter[1] and Hasanjan Sayit

*Cornell University, Cornell University and Worcester Polytechnic Institute*

We show that with suitable restrictions on allowable trading strategies, one has no arbitrage in settings where the traditional theory would admit arbitrage possibilities. In particular, price processes that are not semimartingales are possible in our setting, for example, fractional Brownian motion.

**1. Introduction.** In 1997, Rogers [18] showed that fractional Brownian motion could not be used as a price process for a risky security without introducing arbitrage opportunities. In related work, Delbaen and Schachermayer in 1994 ([5], see also [7]; also [19] and [11] present expository treatments) clarified the concept of no arbitrage by introducing their fundamental concept of "No Free Lunch With Vanishing Risk" (NFLVR) and inter alia showed that as a consequence of the Bichteler–Dellacherie theorem, a necessary condition for a price process to have NFLVR was that it be a semimartingale (see [19] for a nice exposition of this). This insight clarifies the situation of fractional Brownian motion illustrated by Rogers, since fractional Brownian motion is not a semimartingale for most parameter values. Subsequently in his thesis, Cheridito [4] showed that if one properly restricts the class of permissible trading strategies, one can use fractional Brownian motion as a price process and still maintain NFLVR. To accomplish this task, his restriction effectively eliminates those strategies that Rogers had used to illustrate arbitrage. Continuing this line of inquiry, when restricting trading strategies in a manner similar to that proposed by Cheridito, this article attempts to find a general class of processes, which need not be semimartingales that do not permit arbitrage.

The idea is to disallow continuous trading, and moreover to require a minimal fixed time between successive trades. The fixed time can be as

Received June 2008.
[1]Supported in part by NSF Grant DMS-02-02958 and NSA Grant H98230-06-1-0079.

small as one likes, but once chosen, it cannot be changed. This disallows a clustering of trades around a fleeting arbitrage opportunity, such as might occur from a drift process that the random generating process cannot "see." An example might be Brownian local time and Brownian motion, where since the support of the local time is on a (random) set of Lebesgue measure zero, the Brownian motion cannot see when it changes.

This line of inquiry is important for two reasons. First, price processes which are not semimartingales are appearing more regularly in the empirical literature estimating stock price processes (see Lo [14] and references therein), and our methods would provide tools that can be used to determine whether these more general processes are consistent with NFLVR. Second, recent derivatives research has emphasized the importance of transaction costs and illiquidities on restricting the class of permissible trading strategies (see Soner, Shreve and Cvitanic [20] or Jarrow and Protter [12]), for example, continuous trading strategies generate infinite transaction costs under reasonable models of such costs, and those of unbounded variation generate infinite liquidity costs in any finite time interval. As such, these trading strategies could never be used in practice, even if it were physically possible to trade continuously. Without modeling these trading costs explicitly, restricting the class of trading strategies as done by Cheridito provides us with a market setting that implicitly incorporates these trading costs, but maintains the analytic tractability of frictionless markets. Thus, our paper finds those price processes consistent with NFLVR when possible, but also goes beyond those price processes having NFLVR to encompass a larger class. It does this for an extended class of derivative pricing models, without explicitly incorporating transaction costs and illiquidities.

To state the main result of this paper, we need two quick definitions [and we let $\mathbb{F} = (\mathcal{F}_t)_{t \geq 0}$ denote the underlying filtration satisfying the "usual hypotheses"].

DEFINITION 1. The set of simple predictable integrands with bounded support is given by $S(\mathbb{F}) = \{g_0 1_{\{0\}} + \sum_{j=1}^{n-1} g_j 1_{(\tau_j, \tau_{j+1}]} : n \geq 2, 0 \leq \tau_1 \leq \cdots \leq \tau_n$, where all of the $\tau_j$ are $\mathbb{F}$-stopping times; $g_0$ is a real number, and the $g_j$ are real $\mathcal{F}_{\tau_j}$ measurable random variables and $\tau_n$ is bounded$\}$.

We give the name *Cheridito Class* to the trading strategies defined next; we abbreviate it as the class *CC*.

DEFINITION 2 (Cheridito class of trading strategies). For any $h > 0$, let $S^h(\mathbb{F}) = \{g_0 1_{\{0\}} + \sum_{j=1}^{n-1} g_j 1_{(\tau_j, \tau_{j+1}]} \in S(\mathbb{F}) : \forall j, \tau_{j+1} \geq \tau_j + h\}$ and let $\Pi(\mathbb{F}) = \bigcup_{h>0} S^h(\mathbb{F})$. $\Pi$ is the class *CC* of trading strategies.

Let $K^\Pi = \{(H \cdot S)_\infty | H \in \Pi(\mathbb{F})\}$ denote the outcome of the corresponding trading strategies for the price process $S$.

DEFINITION 3. We say $S$ satisfies *the no arbitrage property with respect to the Cheridito Class* $\Pi(\mathbb{F})$ if $K^\Pi \cap L_0^+ = \{0\}$.

THEOREM 1. *Let $S = (S_t)_{t \geq 0}$ be a continuous semimartingale that satisfies the NFLVR property with respect to general admissible integrands such that*

$$(*) \qquad [S,S]_{t+h} - [S,S]_t \geq \delta(h)$$

*for all $t \geq 0$ and any $h > 0$, and for a positive nonrandom increasing function $\delta(\cdot)$ with $\delta(0) = 0$ and $\delta(h) > 0$ for $h > 0$. Assume $[S,S]_t$ is bounded for each $t$. Then for any adapted càdlàg,[2] process $V$ which is bounded in $[0,T]$ for each finite $T > 0$, the process $Y = S + V$ does not have arbitrage in $\Pi(\mathbb{F})$.*

We remark that there are essentially no hypotheses on the bounded process $V$ other than it be càdlàg and adapted to the underlying filtration; for example, it need not have paths of bounded variation. However, simple examples show that the requirement that $V$ be bounded is key.

In this paper, we will also establish related results, consider a different but still restrictive class of trading strategies and prove some useful tools that will allow us to exploit Theorem 1 and give some important examples.

**2. Theorems.** As stated in the Introduction, we will assume given a complete, filtered probability space $(\Omega, \mathcal{F}, P, \mathbb{F})$ satisfying the "usual hypotheses" (i.e., the filtration of $\sigma$ algebras $\mathbb{F}$ is right continuous, and $\mathcal{F}_0$ contains all of the $P$ null sets of $\mathcal{F}$). Let $\mathcal{M}(P)$ be the collection of probability measures on $(\Omega, \mathcal{F})$ that are equivalent to $P$ and we set $L^0 = L^0(\Omega, \mathcal{F}, P)$, $L_{++}^0 = \{\eta \in L^0 : P(\eta \geq 0) = 1 \text{ and } P(\eta > 0) > 0\}$, $L_{--}^0 = \{\eta \in L^0 : P(\eta \leq 0) = 1 \text{ and } P(\eta < 0) > 0\}$, and $L_{+-}^0 = L^0 \setminus (L_{++}^0 \cup L_{--}^0)$.

We begin with a lemma which is the key tool in our analysis. It gives a necessary and sufficient condition for a process in class $CC$ to have no arbitrage.

LEMMA 1. *A process $X_t, t \in [0, \infty)$ satisfies the no arbitrage property in $\Pi(F)$ if and only if for any two bounded stopping times $\tau_1 \geq \tau_0 + h$ with $h > 0$, and any $A \in \mathcal{F}_{\tau_0}$ we have $1_A(X_{\tau_1} - X_{\tau_0}) \in L_{+-}^0$.*

---

[2]Càdlàg is the French acronym for "right continuous with left limits."

PROOF. Let $P_A(\cdot)$ denote the measure $P_A(\cdot) = \frac{P(\cdot \cap A)}{P(A)}$ for any $A \in \mathcal{F}$ with $P(A) > 0$. Let $\tau_0, \tau_1$ be two stopping times with $\tau_1 \geq \tau_0 + h$. Assume X does not have arbitrage but $P_A(X_{\tau_1} > X_{\tau_0}) = 0$; then $X_{\tau_1} \leq X_{\tau_0}$ a.s. on $A$. If $X_{\tau_1} < X_{\tau_0}$ with positive probability on $A$, we take $V = -1_A 1_{(\tau_0, \tau_1]} \in S(\mathbb{F})$ and it is an arbitrage strategy for $X$. This leaves us the only possibility $X_{\tau_1} = X_{\tau_0}$ a.s. on $A$, which is $P_A(X_{\tau_1} = X_{\tau_0}) = 1$. If $P_A(X_{\tau_1} < X_{\tau_0}) = 0$ in the same way, we can show $P_A(X_{\tau_1} = X_{\tau_0}) = 1$. Now, if $P_A(X_{\tau_1} = X_{\tau_0}) < 1$, then either $P_A(X_{\tau_1} > X_{\tau_0}) > 0$ and $P_A(X_{\tau_1} < X_{\tau_0}) > 0$ or one of them is zero. But if one of them is zero, we should have $P_A(X_{\tau_1} = X_{\tau_0}) = 1$ as we showed above and this contradicts $P_A(X_{\tau_1} = X_{\tau_0}) < 1$. So, if $P_A(X_{\tau_1} = X_{\tau_0}) < 1$, then $P_A(X_{\tau_1} > X_{\tau_0}) > 0$ and $P_A(X_{\tau_1} < X_{\tau_0}) > 0$. This proves the sufficiency.

To prove the necessary part, assume there is $V = g_0 1_{\{0\}} + \sum_{j=1}^{n-1} g_j 1_{(\tau_j, \tau_{j+1}]} \in S(F)$ with $P(g_j \neq 0) > 0$ for some $j \in \{1, 2, \ldots, n-1\}$ such that $(V \cdot X)_T \geq 0$ a.s. and $P((V \cdot X)_T > 0) > 0$. Let

$$k = \min\left\{ l : P(g_l \neq 0) > 0, P\left(\sum_{j=1}^{l} g_j (X_{\tau_{j+1}} - X_{\tau_j}) \geq 0\right) = 1, \right.$$

$$\left. P\left(\sum_{j=1}^{l} g_j (X_{\tau_{j+1}} - X_{\tau_j}) > 0\right) > 0 \right\}.$$

If $k = 1$, then $P(g_1 \neq 0) > 0$ and $g_1(X_{\tau_2} - X_{\tau_1}) \geq 0$ a.s. and $g_1(X_{\tau_2} - X_{\tau_1}) > 0$ with positive probability. Let $C = \{g_1(X_{\tau_2} - X_{\tau_1}) > 0\}$ and $A_1 = \{g_1 > 0\}, A_2 = \{g_1 < 0\}$. Then $P(C) > 0$ and either $P(C \cap A_1) > 0$ or $P(C \cap A_2) > 0$. So, assume $P(C \cap A_1) > 0$, since $A_1 \in \mathcal{F}_{\tau_1}$ by hypothesis either $X_{\tau_2} = X_{\tau_1}$ a.s. on $A_1$ or $P_{A_1}(X_{\tau_2} > X_{\tau_1}) > 0$ and $P_{A_1}(X_{\tau_2} < X_{\tau_1}) > 0$. But since $P(A_1 \cap C) > 0$, both $X_{\tau_2} = X_{\tau_1}$ a.s. on $A_1$ and $P_{A_1}(X_{\tau_2} < X_{\tau_1}) > 0$ cannot happen, this contradicts with the hypothesis. If we assume $P(C \cap A_2) > 0$, we also reach the same contradiction, so $k > 1$.

Now, if $k > 1$ then either

$$\sum_{j=1}^{k-1} g_j (X_{\tau_{j+1}} - X_{\tau_j}) \leq 0 \quad \text{a.s.} \quad \text{or} \quad \sum_{j=1}^{k-1} g_j (X_{\tau_{j+1}} - X_{\tau_j}) < 0 \quad \text{a.s.}$$

with positive probability. First assume $\sum_{j=1}^{k-1} g_j (X_{\tau_{j+1}} - X_{\tau_j}) \leq 0$ a.s. and let $A_1 = \{g_k > 0\}$ and $A_2 = \{g_k < 0\}$. Since $P(g_k \neq 0) > 0$, we have either $P(A_1) > 0$ or $P(A_2) > 0$. If $P(A_1) > 0$ and $P(A_2) = 0$, then $X_{\tau_1} = X_{\tau_2}$ a.s. on $A_1$ cannot happen because $P(\sum_{j=1}^{k} g_j (X_{\tau_{j+1}} - X_{\tau_j}) > 0) > 0$ so $P_{A_1}(X_{\tau_{k+1}} > X_{\tau_k}) > 0$ and $P_{A_1}(X_{\tau_{k+1}} < X_{\tau_k}) > 0$ and this contradicts $P(\sum_{j=1}^{k} g_j (X_{\tau_{j+1}} - X_{\tau_j}) \geq 0) = 1$. If $P(A_1) = 0$ and $P(A_2) > 0$ by the same argument as above, we can find a contradiction. If $P(A_1) > 0$ and $P(A_2) > 0$, then both $P_{A_1}(X_{\tau_{k+1}} = X_{\tau_k}) = 1$ and $P_{A_2}(X_{\tau_{k+1}} = X_{\tau_k}) = 1$ cannot happen at the same time. So, either $p_{A_1}(X_{\tau_{k+1}} > X_{\tau_k}) > 0$ and $P_{A_1}(X_{\tau_{k+1}} < X_{\tau_k}) > 0$ or $P_{A_2}(X_{\tau_{k+1}} > X_{\tau_k}) > 0$ and $P_{A_2}(X_{\tau_{k+1}} < X_{\tau_k}) > 0$. In either case, it contradicts $P(\sum_{j=1}^{k} g_j(X_{\tau_{j+1}} - X_{\tau_j}) \geq 0) = 1$.

Now assume $\sum_{j=1}^{k-1} g_j(X_{\tau_{j+1}} - X_{\tau_j}) < 0$ with positive probability. Let $C = \{\sum_{j=1}^{k-1} g_j(X_{\tau_{j+1}} - X_{\tau_j}) < 0\}$, then $P(C) > 0$. And we have $P_C(g_k \neq 0) > 0$ because if $g_k = 0$ a.s. on $C$, then $P(\sum_{j=1}^{k} g_j(X_{\tau_{j+1}} - X_{\tau_j}) < 0) > 0$ a contradiction. So, we either have $P(C \cap A_1) > 0$ or have $P(C \cap A_2) > 0$. We know both $C \cap A_1$ and $C \cap A_2$ are in $\mathcal{F}_{\tau_k}$. Now if $P(C \cap A_1) > 0$, then $P_{C \cap A_1}(X_{\tau_{k+1}} = X_{\tau_k}) = 1$ cannot happen because if not we will have $\sum_{j=1}^{k} g_j(X_{\tau_{j+1}} - X_{\tau_j}) < 0$ on $C \cap A_1$, a contradiction. So, we have $P_{C \cap A_1}(X_{\tau_{k+1}} > X_{\tau_k}) > 0$ and $P_{C \cap A_1}(X_{\tau_{k+1}} < X_{\tau_k}) > 0$ by hypothesis and this contradicts $P(\sum_{j=1}^{k} g_j(X_{\tau_{j+1}} - X_{\tau_j}) \geq 0) = 1$. If $P(C \cap A_2) > 0$ analogously, we reach a contradiction and complete the proof. □

This lemma allows us to prove that the collection of price processes without arbitrage for the class $CC$ is *closed under composition with strictly monotonic functions*.

THEOREM 2. *Let $S = (S_t)_{t \geq 0}$ be a càdlàg stochastic process adapted to the filtration $\mathbb{F}$ and let $f$ be any strictly increasing or strictly decreasing real-valued function in a domain of the real line that contains the range of $X$. Then the no arbitrage property of $X$ in the class $CC$ or $\Pi(\mathbb{F})$ is equivalent to the no arbitrage property of $Y_t = f(S_t)$ in the class $CC$ or $\Pi(\mathbb{F})$, under any measure $Q \in \mathcal{M}(P)$.*

PROOF. Since $X_t$ satisfies the no arbitrage property by Lemma 1 for any $\tau_1 \geq \tau_0$ and any $A \in \mathcal{F}_{\tau_0}$, we have $1_A(X_{\tau_1} - X_{\tau_0}) \in L^0_{+-}(\Omega, \mathcal{F}, P)$. This implies that $1_A(f(X_{\tau_1}) - f(X_{\tau_0})) \in L^0_{+-}(\Omega, \mathcal{F}, P)$ for any strictly monotone function $f$. So, again by the above lemma, the process $f(X_t)$ also satisfies the no arbitrage property. Since $X_t = g(f(X_t))$ for the inverse function $g$ of $f$, by the same argument we know the no arbitrage property of $f(X_t)$ also implies the no arbitrage property of $X_t$. Also, we have $\xi \in L^0_{+-}(\Omega, \mathcal{F}, P)$ if and only if $\xi \in L^0_{+-}(\Omega, \mathcal{F}, Q)$ when $Q$ is equivalent to $P$. So, the claim in Theorem 2 is true for any $Q$ equivalent to $P$. □

Since Cheridito (see [4]) has shown that fractional Brownian motion has no arbitrage in class $CC$, we have the following corollary.

COROLLARY 1. *Let $f$ be any strictly increasing or decreasing function on $\mathbb{R}$ and $B_t^H$ be fractional Brownian motion with Hurst parameter $H$. Then the process $X_t = f(B_t^H)$ does not have arbitrage in class $CC$, that is, $\Pi(\mathbb{F})$ where $\mathbb{F}$ is the natural filtration of $B_t^H$.*

EXAMPLE 1. Let $f(x) = e^x$, then obviously $f$ is a strictly increasing function on $\mathbb{R}$. By Theorem 1, the geometric fractional Brownian motion process $e^{B_t^H}, 0 \leq t \leq T$ does not have arbitrage in class $CC$.

The next theorem shows that the property of no arbitrage for class $CC$ is preserved under filtration shrinkage. Suppose $\mathbb{D}$ is another filtration satisfying the usual hypotheses, and that $\mathcal{D}_t \subset \mathcal{F}_t$ for every $t \geq 0$. We have the following theorem.

THEOREM 3. *Let $X$ be a continuous process adapted to the filtration $\mathbb{D}$, and hence also to $\mathbb{F}$. If $X$ satisfies the no arbitrage property in $\Pi(\mathbb{F})$, then it also satisfies the no arbitrage property in $\Pi(\mathbb{D})$ as well.*

PROOF. Take any two bounded stopping times $\tau_1 \geq \tau_0 + h, h > 0$, of the filtration $\mathbb{D}$. Since $\mathbb{D}$ is the subfiltration of $\mathbb{F}$, $\tau_0$ and $\tau_1$ are also bounded stopping times with respect to $\mathbb{F}$. Then the no arbitrage property of $X$ in $\Pi(\mathbb{F})$ implies, for any $A \in \mathcal{F}_{\tau_0}$, by Lemma 1, we have $1_A(X_{\tau_1} - X_{\tau_0}) \in L^0_{+-}$. Since $\mathbb{D}_{\tau_0} \subset \mathcal{F}_{\tau_0}$ using again Lemma 1, the result is established.

Before restating the main theorem announced in the Introduction and proving it, we establish a second lemma which we will use in its proof. Note that the hypothesis that the stopping times be bounded is essential to the truth of the next lemma. □

LEMMA 2. *Let $B = (B_t)_{t \geq 0}$ be a Brownian motion. Let $T$ denote a finite horizon time and let $0 \leq \tau_1 \leq T$ and $0 \leq \tau_2 \leq T$ be two stopping times with $\tau_2 \geq \tau_1 + h$ for some $h > 0$. Then for any $C > 0$ and any $A \in \mathcal{F}_{\tau_1}$ with $P(A) > 0$, we have $P_A(B_{\tau_2} - B_{\tau_1} > C) > 0$ and $P_A(B_{\tau_2} - B_{\tau_1} < -C) > 0$.*

PROOF. First, we prove $P(\sup_{s \in [h,T]} B_s < -C) > 0$ and $P(\inf_{s \in [h,T]} B_s > C) > 0$. By symmetry, it is enough to prove one of these. So, to show $P(\sup_{s \in [h,T]} B_s < -C) > 0$ note that $\{B_h < -2C, \sup_{s \in [h,T]}(B_s - B_h) < C\} \subseteq \{\sup_{s \in [h,T]} B_s < -C\}$. Using the independence of the increments of Brownian motion, and the fact that $B_h$ is Gaussian and has support over all of $\mathbb{R}$, this implies

$$P\left(\sup_{s \in [h,T]} B_s < -C\right) \geq P\left(B_h < -2C, \sup_{s \in [h,T]}(B_s - B_h) < C\right)$$

$$= P(B_h < -2C) P\left(\sup_{s \in [h,T]}(B_s - B_h) < C\right)$$

$$= P(B_h < -2C) P\left(\sup_{s \in [0,T-h]} B_s < C\right) > 0.$$

Let $Y_t = B_{\tau_1+t} - B_{\tau_1}$. Since Brownian motion "starts afresh" at stopping times, $Y_t$ itself is also a standard Brownian motion and it is independent from $\mathcal{F}_{\tau_1}$. So, $P_A(\sup_{t\in[h,T]} Y_t < -C) = P(\sup_{t\in[h,T]} Y_t < -C) > 0 =$ as shown above (recall that the event $A$ is independent from $Y$). But since $\tau_2 - \tau_1 = \nu$ is a positive random variable and $\nu \geq h$, we have that its values are in the interval $[h, T]$, where $P_A(Y_{\tau_2-\tau_1} < -C) = P_A(Y_\nu < -C) > 0$. And this is equivalent to $P_A(B_{\tau_2} - B_{\tau_1} < -C) > 0$ as required. The other part can be proved by analogously. □

We now prove Theorem 1 which is stated in the Introduction.

PROOF OF THEOREM 1. We first observe that since $S$ is a continuous semimartingale satisfying NFLVR, we can change to an equivalent probability measure such that $S$ is a $\sigma$ martingale; however, since $S$ has continuous paths, we can assume it is, in fact, a local martingale. Therefore, we only need to prove the following:

*Let $X$ be an adapted process on $[0,\infty)$. Assume $X$ has a decomposition $X = M + V$ where $M$ is a continuous local martingale and $V$ is any adapted process. Further assume $M$ satisfies $(*)$, with both $[M, M]$ and $V$ bounded a.s. on $[0, t]$, for each $t \geq 0$. Then $X$ does not have arbitrage in $\Pi(\mathbb{F})$.*

Following the idea of Lemma 2, we let $\tau_1$ and $\tau_2$ be two stopping times with $\tau_2 \geq \tau_1 + h$, and $0 \leq \tau_1 < \tau_2 \leq T$. We only need to show that for any event $A \in \mathcal{F}_{\tau_1}$ with $P(A) > 0$ we have $P_A(X_{\tau_2} > X_{\tau_1}) > 0$ and $P_A(X_{\tau_2} < X_{\tau_1}) > 0$. Since $V$ is a bounded process, the above is satisfied if we can show $P_A(M_{\tau_2} - M_{\tau_1} > C) > 0$ and $P_A(M_{\tau_2} - M_{\tau_1} < -C) > 0$ for any $C >$. Let $\eta_s = \inf\{t > 0 : [M,M]_t > s\}$; for each $s$, $\eta_s$ is a stopping time for $\mathbb{F}$. Condition $(*)$ implies $[M,M]_t \to \infty$ when $t \to \infty$. So, if we let $B_s = M_{\eta_s}$ and $\beta_s = \mathcal{F}_{\eta_s}$, then $(B_s, \beta_s)$ is a standard Brownian motion, and moreover, $M_t = B_{[M,M]_t}$ (cf., e.g., [16]). Since $[M, M]$ has continuous and strictly increasing paths, we have $\eta_{[M,M]_t} = t$, and of course $\{[M,M]_u \leq s\} = \{\eta_s \geq u\} \in \mathcal{F}_{\eta_s}$ and both $[M,M]_{\tau_1}$ and $[M,M]_{\tau_2}$ are stopping times for the filtration $\beta$. Last, note that $\mathcal{F}_{\tau_1} \subseteq \beta_{[M,M]_{\tau_1}}$, so if $A$ is an event in $\mathcal{F}_{\tau_1}$, then also $A \in \beta_{[M,M]_{\tau_1}}$. Since $[M, M]$ is bounded, so too are the stopping times $[M,M]_{\tau_1}$ and $[M,M]_{\tau_2}$ and, therefore, by the hypothesis $(*)$ we have $[M,M]_{\tau_2} - [M,M]_{\tau_1} \geq \delta(\tau_2 - \tau_1) \geq \delta(h)$. Then by Lemma 2 we have for any $C > 0$

$$P_A(M_{\tau_2} - M_{\tau_1} \geq C) = P_A(B_{[M,M]_{\tau_2}} - B_{[M,M]_{\tau_1}} \geq C) > 0$$

and

$$P_A(M_{\tau_2} - M_{\tau_1} \leq -C) = P_A(B_{[M,M]_{\tau_2}} - B_{[M,M]_{\tau_1}} \leq -C) > 0$$

and the theorem is proved. □

**3. Examples.**

EXAMPLE 2. Let $X$ be given by $X_t = \int_0^t B_s\, dB_s + t$ for $t \geq 0$. In this example, the quadratic variation process is $\int_0^t B_s^2\, ds$ which is strictly increasing a.s., and the process $V_t = t$ is a bounded process on $[0, t]$ for each $t$. However, by Itô's formula, we have $X_t = \frac{1}{2}(B_t^2 + t)$, which has arbitrage in the Cheridito class $\Pi(\mathbb{F})$, and in any other reasonable framework. Here, the martingale term $M_t = \int_0^t B_s\, dB_s$ does not satisfy $(*)$. This shows that the condition $(*)$ cannot be easily improved upon.

EXAMPLE 3. Let $X_t = \int_0^t \mu_s\, dB_s + V_t$ be a continuous process with $P(\int_0^T \mu_s^2\, ds < \infty) = 1$ and $\inf_{s \in [0,T]} |\mu_s| \geq \delta$ for some $\delta > 0$ and assume $V_t$ and $\mu_t$ are bounded adapted processes. Then X satisfies the no arbitrage property in $\Pi(\mathbb{F})$. This follows because $M_t = \int_0^t \mu_s\, dB_s$ and $[M, M]$ satisfies $(*)$.

EXAMPLE 4. By Tanaka's formula, we have $|B_t| = \int_0^t \text{sign}(B_s)\, dB_s + L_t$, where $L_t$ is the local time at zero of the Brownian motion. Since $|B_t|$ is a positive process beginning from zero, it has arbitrage in $\Pi$. The local martingale part $\int_0^t \text{sign}(B_s)\, dB_s$ satisfies condition $(*)$, but $L_t$ is an unbounded process in $[0, T]$.

EXAMPLE 5. In Example 4, we now let $\tau = \inf\{t > 0 | L_t > N\}$ for positive $N$ and let $D_t = \int_0^t \text{sign}(B_s)\, dB_s + L_{t \wedge \tau}$. Then this process, modified from example 4, does not have arbitrage in $\Pi$.

EXAMPLE 6. Consider the processes $X_t^\alpha = \int_0^t s^\alpha\, dB_s + V_t$ on $[0, T]$. Here, $\alpha > -\frac{1}{2}$ and $V_t$ is any adapted bounded process with respect to Brownian motion. The quadratic variation of $M_t = \int_0^t s^\alpha\, dB_s$ is $[M, M]_t = \int_0^t s^{2\alpha}\, ds$. By a simple calculation, we get $[M, M]_{t+h} - [M, M]_t = \frac{1}{2\alpha+1}[(t+h)^{2\alpha+1} - t^{2\alpha+1}] \geq \frac{1}{2\alpha+1} h^{2\alpha+1}$. So, we can let $\delta(h) = \frac{1}{2\alpha+1} h^{2\alpha+1}$ and the condition of the theorem is satisfied. So, the processes $X_t^\alpha$ do not have arbitrage on $\Pi(\mathbb{B})$.

The remainder of this section will be devoted to the examples within the class of *Gaussian moving average processes*, which will include the case of *fractional Brownian motions*. This treatment will allow us, en passant, to give a new proof of Cheridito's result that fractional Brownian motion does not allow arbitrage in $CC$ [4]. What underlies this treatment is the theorem that the Delbaen–Schachermayer condition on the price process of NFLVR implies that the price process must be a semimartingale.

We consider a probability space equipped with a two-sided Brownian motion $(W_t)_{t \in \mathbb{R}}$, that is, $W$ is a continuous centered Gaussian process with covariance

$$\operatorname{Cov}(W_t, W_s) = \tfrac{1}{2}(|t| + |s| - |t - s|), \qquad t, s \in \mathbb{R}.$$

For any function $\varphi : \mathbb{R} \to \mathbb{R}$ that is zero on the negative real axis and satisfies for all $t > 0$,

(1)
$$\varphi(t - \cdot) - \varphi(-\cdot) \in L^2(R, R),$$
$$Y_t^\varphi = \int_{-\infty}^t [\varphi(t - u) - \varphi(-u)] \, dW_u, \qquad t \in \mathbb{R}.$$

We recall that a stochastic process $(Y_t)_{t \in \mathbb{R}}$ has stationary increments if for all $t_0 \in \mathbb{R}$,

$$(Y_{t+t_0} - Y_{t_0})_{t \in R} \stackrel{d}{=} (Y_t - Y_{t_0})_{t \in \mathbb{R}}.$$

Cheridito showed in [4] that the process $Y_t^\varphi$ is a semimartingale in $[0, T]$ for some $T > 0$ if and only if $\varphi$ has the following form:

(2)
$$\varphi(t) = \begin{cases} v + \int_0^t \psi(s) \, ds, & \text{when } t \geq 0, \\ 0, & \text{when } t < 0, \end{cases}$$

where $\psi \in L^2(\mathbb{R}_+, \mathbb{R})$ and $v \in \mathbb{R}$. A key example is that if we let $\varphi(t) = 1_{(0,\infty)} t^{H-1/2}, t \in R$ for $H \in (0,1)$, then $\varphi(t)$ does not satisfy equation (2), where $Y^\varphi$ is not a semimartingale. We note that

$$\operatorname{Cov}(Y_t^\varphi, Y_s^\varphi) = \tfrac{1}{2} c_H^2 (|t|^{2H} + |s|^{2H} - |t - s|^{2H}), \qquad t, s \in \mathbb{R},$$

where

$$c_H = \left( \frac{1}{2H} + \int_0^\infty [(1+u)^{H-1/2} - u^{H-1/2}] \, du \right)^{1/2}.$$

These processes are called *fractional Brownian motions*. Since these processes are not semimartingales, they cannot satisfy NFLVR, hence by the definition of NFLVR, we can conclude that there must exist a sequence $H_n$ of simple predictable processes of bounded support such that

$$(H_n \cdot Y^\varphi)_\infty \geq -\tfrac{1}{2}, \qquad (H_n \cdot S)_\infty \to f$$

for a function $f \geq 0, P(f > 0) > 0$. (The choice of $\tfrac{1}{2}$ is, of course, arbitrary.) We recall the Dalang–Morton–Willinger theorem.

THEOREM 4. *Let $(S_n)_{n=0,1,\ldots,N}$ be a process adapted to $(\mathcal{F}_n)_{n=0,1,\ldots,N}$. Let*

$$H = \{\sum_{n=1}^N f_{n-1}(S_n - S_{n-1}),$$

where each $f_n : \Omega \to \mathbb{R}$ is $\mathcal{F}_n$ measurable. If $H \cap L^0_+ = \{0\}$, then there is an equivalent measure $Q$ such that $S$ is a $Q$-martingale.

An easy consequence of this theorem in our setting is the following corollary.

COROLLARY 2. *A process $X$ satisfies the no arbitrage property in $\Pi(\mathbb{F})$ if and only if for any sequence of bounded stopping times that satisfies $\tau_1 \leq \tau_2 \leq \cdots \leq \tau_N; \tau_{i+1} \geq \tau_i + h, i = 1, 2, \ldots, N - 1$, for some $h > 0$, there is an equivalent probability measure $Q$ such that $(X_i, \mathcal{F}_i), i = 1, 2, \ldots, N$, is a $Q$-martingale.*

We will need the following elementary lemma which we found in [15].

LEMMA 3. *Let $(S_i, \mathcal{F}_i)_{\{i=0,1,2,\ldots,N\}}$ be an adapted real valued process such that for every predictable process $(h_i)_{\{i=0,1,2,\ldots,N\}}$ we have that $(h \cdot S)_N = \sum_{i=1}^N h_i \triangle S_i$ is unbounded from above and from below as soon as $(h \cdot S)_N \neq 0$. Then there is a measure $Q$ equivalent to the original measure such that $(S_i)_{\{i=0,1,2,\ldots,N\}}$ is a $Q$-martingale for the underlying filtration $\mathbb{F} = \mathcal{F}_{\{i=0,1,2,\ldots,N\}}$.*

The key lemma for this topic is as follows.

LEMMA 4. *Let $X = (X_t)_{t \geq 0}$ with filtration $\mathbb{F}$ be an adapted continuous process and $\tau$ be any bounded stopping time. If for any $A \in \mathcal{F}_\tau$ with $P(A) > 0$ and any $0 < \delta < T < \infty$, we have $P(\{1_A \sup_{t \in [\delta, T]}(X_{\tau+t} - X_\tau) < -C\}) > 0$ and $P(\{1_A \inf_{t \in [\delta, T]}(X_{\tau+t} - X_\tau) > C\}) > 0$ is satisfied for all $C > 0$, then for any bounded adapted process $V$, $Y = X + V$ does not have arbitrage in $\Pi(\mathbb{F})$.*

PROOF. Fix any sequence of bounded stopping times $\tau_1 \leq \tau_2 \leq \cdots \leq \tau_{N+1}, \tau_{i+1} \geq \tau_i + h, i = 1, 2, \ldots, N$, for some $h > 0$. By Corollary 2, we need to show there is an equivalent probability measure $Q$ such that $(Y_{\tau_i}, \mathcal{F}_{\tau_i})_{i=1}^{N+1}$ is a martingale under $Q$. We prove this using Lemma 3. So take any nontrivial predictable simple process

$$H = \sum_{i=1}^{N} g_i 1_{(\tau_i, \tau_{i+1}]}.$$

We assume $g_n \neq 0$. Consequently, either $P(g_n > 0) > 0$ or $P(g_n < 0) > 0$. So, we assume $P(g_n > 0) > 0$. We can choose a big enough number $M > 0$ such that the event $A = (\{\sum_{i=1}^{n-1} g_i(Y_{\tau_{i+1}} - Y_{\tau_i}) < M\} \cap \{g_n > 0\}$ has positive probability, namely $P(A) > 0$. We note that $A \in \mathcal{F}_{\tau_n}$. Then by the hypotheses of this lemma, we have $P(\{1_A \sup_{t \in [1/2h, d]}(X_{\tau+t} - X_\tau) < -C\}) > 0$ and $P(\{1_A \inf_{t \in [1/2h,d]}(X_{\tau+t} - X_\tau) > C\}) > 0$ for any $C > 0$. Here, $d$ is a number greater than the bound of $\tau_{n+1}$. Since $V$ is bounded and $\tau_{n+1} \geq \tau_n + h$, we have $P(1_A(Y_{\tau_2} - Y_{\tau_1}) < -C) > 0$ and $P(1_A(Y_{\tau_2} - Y_{\tau_1}) > C) > 0$ for any $C > 0$. But the sum $\sum_{i=1}^{n-1} g_i(Y_{\tau_{i+1}} - Y_{\tau_i})$ is bounded on $A$, so we have that $(H \cdot Y)_n$ is unbounded from below and above. Then by Lemma 3, there is an equivalent measure $Q$ that makes $Y_{\tau_i}, i = 1, 2, \ldots, N+1$, a martingale. □

THEOREM 5. *Let $Y^\varphi$ be a moving average process as given in equation (1). Let the stationary centered Gaussian process $X_t^\varphi = \int_\mathbb{R} \varphi(t-u) \, dW_u$ satisfy the following property: For any $0 < \delta < T < \infty$, $P(\sup_{t \in [\delta,T]} X_t^\varphi < -C) > 0$ and $P(\inf_{t \in [\delta,T]} X_t^\varphi > C) > 0$ for all $C > 0$. Then, for any bounded process $V$, the process $Z_t = Y_t^\varphi + V_t$ does not have arbitrage in $\Pi(F)$. In particular, the process $Z_t = B_t^H + V_t$ does not have arbitrage in $\Pi(F)$. Here, $B_t^H$ is fractional Brownian motion.*

REMARK 1. We remark that in Lemma 4.2 of [4] Cheridito has shown that when $\varphi = 1_{(0,\infty]}(t) t^{H-1/2}$ for $H \in (0, \frac{1}{2}) \cup (\frac{1}{2}, 1)$, $X_t^\varphi$ satisfies the conditions of the above theorem. Then by using this lemma, he proved that the process $B_t^H + v(t)$, where $v(t)$ is a deterministic function, does not have arbitrage in $\Pi(\mathbb{F})$. Our proof of the result might be considered more simple, and also it extends the result to the random case.

PROOF OF THEOREM 5. Let $\tau$ be any bounded stopping time. By Lemma 4, we only need to check the process $Z_t = 1_A(Y_{\tau+t}^\varphi - Y_\tau^\varphi)$ satisfies $P(\sup_{t \in [\delta,T]} Z_t < -C) > 0$ and $P(\inf_{t \in [\delta,T]} Z_t < -C) > 0$ for any $A \in \Im_\tau$ with $P(A) > 0$ and any $0 < \delta < T < \infty, C > 0$. To prove these, we borrow the idea of the proof of Theorem 4.3 of [4]. Let

$$\tilde{\Omega} := \left\{ \omega \in C(\mathbb{R}) : \omega(0) = 0 \text{ and } \forall t \in \mathbb{R}, \lim_{s \to t} \frac{\omega(t) - \omega(s)}{\sqrt{|t-s| \log(1/|t-s|)}} = 0 \right\}.$$

Let $\mathcal{B}$ be the $\sigma$-algebra of subsets of $\tilde{\Omega}$ that is generated by the cylinder sets, and $P$ be Wiener measure on $(\tilde{\Omega}, \mathcal{B})$. Without loss of generality, we assume that $(Y_t^\varphi)$ is defined on $(\tilde{\Omega}, \mathcal{B}, P)$ by the improper Riemann–Stieltjes integrals

$$Y_t^\varphi(\omega) = \int_{-\infty}^t [\varphi(t-s) - \varphi(-s)] \, d\omega(s).$$

We define the filtration $\mathcal{F}^{\tilde{\Omega}} = (\mathcal{F}_t^{\tilde{\Omega}}), t \in \mathbb{R}$ by

$$\mathcal{F}_t^{\tilde{\Omega}} := \sigma\{\{\omega \in \tilde{\Omega} : \omega(s) \leq a\} : -\infty < s \leq t, a \in \mathbb{R}\}.$$

It is clear that $\mathcal{F}^{\tilde{\Omega}}$ contains the filtration $\mathcal{F}^{Y^\varphi} = (\mathcal{F}_t^{Y^\varphi})_{t \in \mathbb{R}}$, which is given by

$$\mathcal{F}_t^{Y^\varphi} = \sigma\{Y_s^\varphi : 0 \leq s \leq t\}.$$

Therefore, $\tau$ is also $\mathcal{F}^{\tilde{\Omega}}$ stopping time. Now we split each function $\omega \in \tilde{\Omega}$ at the time point $\tau(\omega)$. Let

$$\pi_1\omega(s) := \omega(s)1_{(-\infty,\tau(\omega)]}(s), \qquad s \in \mathbb{R},$$

$$\pi_2\omega(s) := \omega(\tau(\omega) + s) - \omega(\tau(\omega)), \qquad s \geq 0$$

and let

$$\Omega_1 = \{\pi_1(\omega) \in R^R : \omega \in \tilde{\Omega}\},$$

$\mathcal{B}_1$ be the $\sigma$-algebra of subsets of $\Omega_1$ that is generated by the cylinder sets.

$$\Omega_2 = \{\pi_2(\omega) \in C[0,\infty) : \omega \in \tilde{\Omega}\}$$

and $\mathcal{B}_2$ the $\sigma$-algebra of subsets of $\Omega_2$ that is generated by the cylinder sets. It can be easily checked that the mapping $\pi_1 : (\tilde{\Omega}, \mathcal{B}) \to (\Omega_1, \mathcal{B}_\infty)$ is $\mathcal{F}_\tau^{\tilde{\Omega}}$ measurable. On the other hand, since a Lévy process "renews itself at stopping times" (see, e.g., [16], page 23), it follows that $(\pi_2\omega(s))_{s\geq 0}$ is a Brownian motion which is independent of $\mathcal{F}_\tau^{\tilde{\Omega}}$. We have

$$1_A(Y_{\tau+t}^\varphi - Y_\tau^\varphi) = \int_{-\infty}^\tau [\varphi(\tau + t - s) - \varphi(-s)]\,dW_s + \int_0^t \varphi(t-s)\,dW_s.$$

Let

$$U_t(\omega_1, \omega_2) := 1_A(\omega_1) \int_{-\infty}^{\tau(\omega_1)} [\varphi(\tau(\omega_1) + t - s) - \varphi(-s)]\,d\omega_1(s)$$

$$+ 1_A(\omega_1) \int_0^t \varphi(t-s)\,d\omega_2(s)$$

for $\omega_1 \in \Omega_1, \omega_2 \in \Omega_2$ and $t \geq 0$. Then for all $\omega \in \tilde{\Omega}$ and $t \geq 0$, we have the following relation

$$[1_A(Y_{\tau+t}^\varphi - Y_\tau^\varphi)](\omega) = U_t(\pi_1\omega, \pi_2\omega).$$

For each fixed $\omega_1$, the process $1_A(\omega_1)\int_{-\infty}^{\tau(\omega_1)}[\varphi(\tau(\omega_1)+t-s)-\varphi(-s)]\,d\omega_1(s)$ is a continuous process and so

$$\sup_{t \in [\delta,T]} \left(1_A(\omega_1) \int_{-\infty}^{\tau(\omega_1)} [\varphi(\tau(\omega_1) + t - s) - \varphi(-s)]\,d\omega_1(s)\right)$$

is finite. Since $(U_t)_{t\in[[\delta,T]}$ is a continuous stochastic process on $(\Omega_1 \times \Omega_2, \mathcal{B}_1 \times \mathcal{B}_2)$, the set

$$F := \left\{(\omega_1, \omega_2) \in \Omega_1 \times \Omega_2 : \sup_{t \in [\delta,T]} U_t(\omega_1, \omega_2) \leq -C\right\}$$

is $\mathcal{B}_1 \times \mathcal{B}_2$-measurable. It follows (see, e.g., Proposition A.2.5 of [13]) that for almost every $\omega \in \tilde{\Omega}$,

$$E[1_F(\pi_1, \pi_2) \mid \mathcal{F}_\tau^{\tilde{\Omega}}](\omega) = \phi(\pi_1 \omega),$$

where the mapping $\phi: \Omega_1 \to \mathbb{R}$ is given by $\phi(\omega_1) := E[1_A(\omega_1, \pi_2)], \phi_1 \in \Omega_1$. Since $(\pi_2 \omega(t))_{t \geq 0}$ is a Brownian motion under $P$, by the condition of the theorem, it follows that for any fixed $\omega_1 \in A$

$$\phi(\omega_1) = P\left[\sup_{t \in [\delta, T]} U_t(\omega_1, \pi_2) \leq -C\right]$$

$$(3) \quad \geq P\left[\sup_{t \in [\delta, T]} \left(1_A(\omega_1) \int_{-\infty}^{\tau(\omega_1)} [\varphi(\tau(\omega_1) + t - s) - \varphi(-s)] \, d\omega_1(s)\right)\right.$$

$$\left. + 1_A(\omega_1) \sup_{t \in [\delta, T]} \int_0^t \varphi(t - s) \, d\omega_2(s) \leq -C\right] > 0,$$

Therefore,

$$P\left[\sup_{t \in [\delta, T]} 1_A(Y_{\tau+t}^\varphi - Y_\tau^\varphi) \leq -C\right] = E[1_F(\pi_1, \pi_2)] = E[E[1_F(\pi_1, \pi_2) \mid \mathcal{F}_\tau^{\tilde{\Omega}}]]$$

$$(4) \qquad = E[\phi \circ \pi_1] > 0.$$

By using the same argument above, one can prove $P[\inf_{t \in [\delta, T]} 1_A(Y_{\tau+t}^\varphi - Y_\tau^\varphi) \geq C] > 0$. This completes the proof of the theorem. $\square$

As a corollary of the theorem, we state a property of fractional Brownian motion.

COROLLARY 3. *Let $B_t^H, t \in [0, \infty]$ be a fractional Brownian motion. Let $\tau_1$ and $\tau_2$ be any two bounded stopping times with $\tau_2 \geq \tau_1 + h$ for some $h > 0$. Then for any $C > 0$ and any $A \in \mathcal{F}_{\tau_1}$ with $P(A) > 0$, we have $P_A(B_{\tau_2}^H - B_{\tau_1}^H > C) > 0$ and $P_A(B_{\tau_2}^H - B_{\tau_1}^H < -C) > 0$.*

When $\frac{1}{2} < H < 1$, the fractional Brownian motion $B^H$ admits the following integral representation

$$B_t^H = \int_0^t K_H(t, s) \, dB_s,$$

where B is the ordinary Brownian motion and the kernel $K_H(t, s)$ has the form

$$K_H(t, s) = C_H \left[ \left(\frac{t}{s}\right)^{H-1/2} (t-s)^{H-1/2} \right.$$

$$\left. - \left(H - \frac{1}{2}\right) s^{1/2-H} \int_0^t u^{H-3/2} (u-s)^{H-1/2} \, du \right].$$

For each fixed $t > 0$, the kernel $K_H$ defines an operator $\mathcal{K}_H$ in $L([0,t])^2$ by

$$(\mathcal{K}_H h)(u) = \int_0^t K_H(u,s) h(s) \, ds.$$

And for any absolutely continuous $h$, its inverse operator is given by

$$(\mathcal{K}_H^{-1} h)(s) = s^{H-1/2} D_{0+}^{H-1/2} (r^{1/2} h'(r))(s),$$

where $D_{0+}^{H-1/2}$ denotes the left-fractional derivative, defined for $t > 0$ by

$$D_{0+}^{H-1/2} f(t) = \frac{1}{\Gamma(3/2 - H)} \left( \frac{f(t)}{t^{H-1/2}} + \left( H - \frac{1}{2} \right) \int_0^t \frac{f(t) - f(s)}{(t-s)^{H+1/2}} \, ds \right).$$

Then for any adapted process $\mu_t, t \in [0, \infty)$ with $\int_0^t \mu_s \, ds < \infty$ a.s., the process $\tilde{B}_t^H = B_t^H + \int_0^t \mu_s \, ds$ is again a fractional Brownian motion with Hurst parameter $H$ in the interval $[0, T]$, for fixed $T > 0$, under and equivalent change of measure $Q_T$ defined by $\frac{dQ_T}{dp} = \Lambda_T$ if the following two conditions are satisfied

(i) $\int_0^{\cdot} \mu_s \, ds \in$ the image of $\mathcal{K}_H$ a.s.
(ii) $E(\Lambda_t) = 1$ for any $t > 0$, where

$$\Lambda_t = \exp\left( -\int_0^t \left( \mathcal{K}_H^{-1} \int_0^{\cdot} \mu_s \, ds \right)(s) \, dB_s + \frac{1}{2} \int_0^t \left( \mathcal{K}_H^{-1} \int_0^{\cdot} \mu_s \, ds \right)^2(s) \, ds \right).$$

COROLLARY 4. *Assume $\int_0^{\cdot} \mu_s \, ds$ is in the image of $\mathcal{K}_H$, a.s. and $E(\Lambda_t) = 1$ for any $t > 0$. Then the process $X_t^H =: B_t^H + \int_0^t \mu_s \, ds + V_t$ satisfies the no arbitrage property in $\Pi(F)$ when $\frac{1}{2} < H < 1$, where $V$ is any adapted bounded process.*

PROOF. Take any element $H_s = \sum_{i=1}^n C_i 1_{(\tau_i, \tau_{i+1}]}(s)$ in $\Pi(F)$. Since the stopping times $\tau_i, 1 \leq i \leq n+1$ are bounded, the process $H_s$ is supported in some interval $[0, T]$ for a finite $T > 0$. By Girsanov's theorem for fractional Brownian motions, the process $B_t^H + \int_0^t \mu_s \, ds$ is a fractional Brownian motion in the interval $[0, T]$ under the measure $Q_T$. Then by Corollary 3, we have $Q_T((H \cdot X)_T < 0) > 0$. This shows that our claim is true. □

**4. Time changed processes.** Time changes have recently become quite popular in the construction of price processes. The reader can consult, for example, any or all of [2, 9] and [10]. Therefore, it is interesting to check stability of our no arbitrage property under a time change.

In this section, let $\mathbb{F}$ denote the underlying filtration satisfying the usual hypotheses, and let $(\nu_t)_{t \geq 0}$ denote a *change of time*, that is, $(\nu_t)_{t \geq 0}$ is a family of $\mathbb{F}$ stopping times such that $t \to \nu_t(\omega)$ is right continuous and nondecreasing for almost all $\omega$, $\nu_t < \infty$ a.s. and $\nu_0 = 0$. A *continuous change of time* is one where $t \to \nu_t(\omega)$ is continuous for almost all $\omega$. We let $\tilde{\mathbb{F}}$ denote the time changed filtration given by $\tilde{\mathcal{F}}_t = \mathcal{F}_{\nu_t}$.

In order to prove Theorem 6 which follows, we need to remark that Lemma 1 holds under a weaker condition. To be precise, we state the new version. Note that the proof of Lemma 1 can be used to prove Lemma 5 with only slight modifications, and that if a process has no arbitrage for the class of simple predictable integrands of bounded support, then it a fortiori has no arbitrage for the Cheridito Class $(CC)$.

LEMMA 5. *A càdlàg adapted process $(X_t)_{t \geq 0}$ satisfies the no arbitrage property in the class of simple predictable integrands of bounded support if and only if for any bounded stopping times $\tau_1 \geq \tau_0$ and any $A \in \mathcal{F}_{\tau_0}$ we have $1_A(X_{\tau_1} - X_{\tau_0}) \in L^0_{+-}$.*

By using the same reasoning as in the proof of Theorem 2, we can prove the following corollary.

We remark that Delbaen and Schachermayer long ago considered the restriction to simple integrands of bounded support, and showed in 1994 (cf. [5]) that NFLVR for this framework implies the existence of an equivalent local martingale measure (see their Theorem 7.6 of [5]).

COROLLARY 5. *Let $S_t$ be a càdlàg stochastic process adapted to the filtration $\mathbb{F}$ and let $f$ be any strictly increasing or strictly decreasing function in a domain of the real line that contains the range of $S$. Then the no arbitrage property of $S$ in $S(\mathbb{F})$ is equivalent to the no arbitrage property of $Y_t = f(S_t)$ in $S(\mathbb{F})$ under any equivalent change of (probability) measure $Q$.*

Corollary 5 provides a wealth of examples of processes which are not semimartingales (and, therefore, cannot satisfy NFLVR), but which nevertheless have the no arbitrage property in $S(\mathbb{F})$. The following example, taken from [16], is typical.

EXAMPLE 7. Let $B$ be a one dimensional Brownian motion. Then the process $Y_t = B^{1/3}$ satisfies the no arbitrage property in $S(\mathbb{F})$, where $\mathbb{F}$ is the filtration of $Y$ (which is the same as the filtration of the Brownian motion $B$).

Of course, $Y$ is the composition of $B$ with the strictly increasing function $f(x) = x^{1/3}$, and since $B$ satisfies NFLVR, it is clearly a no arbitrage process, and hence so is $Y$ by Theorem 5.

Here, we remark that the critical reason for us to be able to state both Theorem 2 and Corollary 5 is that we allow both short and long positions of any amount in our trading strategies. More precisely for each $H = \sum_i^n H_i 1_{[\tau_i,\tau_{i+1})}$ either in $\Pi(\mathbb{F})$ or in $S(\mathbb{F})$, we allow $H_i$ to be any random variable. If one considers the class of trading strategies that restricts $H_i$ in both $\Pi(\mathbb{F})$ and $S(\mathbb{F})$ to be bounded random variables, one can still get a result as in Theorem 2. But if one adds a short sale restriction, namely if each $H_i$ is only allowed to be a nonnegative random variable, then one can argue that result such as Theorem 2 does not hold in general.

Another consequence of Lemma 5 is the following simple method to check for arbitrage. Due to this lemma, we can replace the class $(CC)$ with the class of simple predictable integrands of bounded support.

COROLLARY 6. *A process $(X_t)_{t\geq 0}$ has arbitrage in the class of simple predictable integrands of bounded support if and only if its arbitrage strategy can be taken in the form $1_{(\tau_0,\tau_1]}$ or $-1_{(\tau_0,\tau_1]}$, for two bounded stopping times $0 \leq \tau_0 \leq \tau_1$.*

PROOF. If the process has arbitrage then by Lemma 5, there are two bounded stopping times $\tau_1 \leq \tau_2$ and a nontrivial event $A \in \mathcal{F}_{\tau_1}$ such that either $P_A(S_{\tau_1} \leq S_{\tau_2}) = 1$ or $P_A(S_{\tau_1} \geq S_{\tau_2}) = 1$ and $P_A(S_{\tau_1} = S_{\tau_2}) < 1$. This implies either $1_A 1_{(\tau_1,\tau_2]}$ or $-1_A 1_{(\tau_1,\tau_2]}$ is an arbitrage strategy. Define stopping times

$$\tau_1^A(\omega) = \begin{cases} M, & \text{when } \omega \notin A, \\ \tau_1, & \text{when } \omega \in A, \end{cases}$$

$$\tau_2^A(\omega) = \begin{cases} M, & \text{when } \omega \notin A, \\ \tau_2, & \text{when } \omega \in A, \end{cases}$$

where $M$ is any number bigger than the bound of $\tau_2$. Then $1_A 1_{(\tau_1,\tau_2]} = 1_{(\tau_1^A,\tau_2^A]}$ and this completes the proof. □

THEOREM 6. *Let $(X_t)_{t\geq 0}$ be an $\mathbb{F}$ adapted process satisfying the no arbitrage property on $S(\mathbb{F})$. Let $(\nu_t)_{t\geq 0}$ be a continuous change of time, such that $\nu_t$ is bounded a.s. for each $t \geq 0$. Let $\tilde{X} = X_{\nu_t}$, and $\tilde{\mathcal{F}}_t = \mathcal{F}_{\nu_t}$, for $t \geq 0$. Then the no arbitrage property of $X$ on $S(F)$ is equivalent to the no arbitrage property of the time changed process $\tilde{X}$ on $S(\tilde{F})$.*

PROOF. By Lemma 5, it suffices to check that for any two bounded stopping times $\tau_0 \leq \tau_1$ of $\tilde{\mathbb{F}}$ and for any $A \in \tilde{\mathcal{F}}_{\tau_0}$ we have either $P_A(\tilde{X}_{\tau_0} > \tilde{X}_{\tau_1})$ and $P_A(\tilde{X}_{\tau_0} < \tilde{X}_{\tau_1})$ or $P_A(\tilde{X}_{\tau_0} = \tilde{X}_{\tau_1}) = 1$. To do this, we first define $C_s = \inf\{t > 0 : \nu_t > s\}$. Since $\nu_t$ is continuous, we have $\nu_{C_s} = s$. Note that all of $C_s$ are $\tilde{\mathbb{F}}$ stopping times, and for any stopping time $\tau$ of $\tilde{\mathbb{F}}$, since $\{\nu_\tau \geq s\} = \{C_s \leq \tau\} \in \tilde{\mathcal{F}}_{C_s} = \mathcal{F}_s$, we know $\nu_\tau$ is an $\mathbb{F}$ stopping time. So, $\nu_{\tau_0}$ and $\nu_{\tau_1}$ are bounded stopping times of $\mathbb{F}$, and we have $\tilde{\mathcal{F}}_{\tau_0} = \mathcal{F}_{\nu_{\tau_0}}$. Since $X$ satisfies the no arbitrage property for $S(\mathbb{F})$, by Lemma 5 we have $P_A(X_{\nu_{\tau_0}} > X_{\nu_{\tau_1}})$ and $P_A(X_{\nu_{\tau_0}} < X_{\nu_{\tau_1}})$ or $P_A(X_{\nu_{\tau_0}} = \tilde{X}_{\nu_{\tau_1}}) = 1$ so the above conditions are satisfied. Now if $X$ has arbitrage on $S(\mathbb{F})$, then by Corollary 6 the arbitrage strategy can be taken in the form $1_{(\tau_0,\tau_1]}$ or $-1_{(\tau_0,\tau_1]}$ for two bounded stopping times $\tau_0 \leq \tau_1$. Then one can easily check that either $1_{(C_{\tau_0},C_{\tau_1}]}$ or $-1_{(C_{\tau_0},C_{\tau_1}]}$ is an arbitrage strategy for $\tilde{X}$ on $S(\tilde{\mathbb{F}})$. $\square$

As an application of Theorem 6, we have the following theorem.

THEOREM 7. *Let $S$ be a semimartingale that admits an equivalent local martingale measure. Assume $[S,S]_t \to \infty$ a.s. when $t \to \infty$ and that $[S,S]_t$ is bounded, each $t \geq 0$. Then*

$$Z_t^\alpha = S_t + [S,S]_t^\alpha$$

*satisfies the no arbitrage property with respect to the class of simple trading strategies of bounded support [and thus a fortiori with respect to the Cheridito Class (CC)] when $\alpha \geq \frac{1}{2}$.*

Before proving Theorem 7, we establish a lemma.

LEMMA 6. *Let $B_t$ be a Brownian motion and $\alpha \geq \frac{1}{2}$. Then the processes $Y_t^\alpha = B_t + t^\alpha$ each satisfy the no arbitrage property with respect to simple predictable processes of bounded support.*

PROOF. Take any two bounded stopping times $0 \leq \tau_0 \leq \tau_1$. By Blumenthal's $0-1$ law, we either have $\tau_0 > 0$ a.s. or $\tau_0 = 0$. If $\tau_0 > 0$, then by Girsanov's theorem, we know the process $Y^\alpha - (Y^\alpha)^{\tau_0}$ admits an equivalent martingale measure $Q$ with density given by $f = \exp(-\alpha \int_{\tau_0}^1 t^{\alpha-1} dB_t - \frac{\alpha^2}{2} \int_{\tau_0}^1 t^{2\alpha-2} dt)$. So, we have $E_Q Y_{\tau_1}^\alpha = E_Q Y_{\tau_0}^\alpha$ and this shows that neither $1_{(\tau_0,\tau_1]}$ nor $-1_{(\tau_0,\tau_1]}$ can be the arbitrage strategy. Therefore, assume $\tau_0 = 0$. By the law of the iterated logarithm, $\liminf_{t \to 0} \frac{B_t}{t^{1/2}\sqrt{2\log\log(1/t)}} = -1$. By writing $Y_t^\alpha = B_t + t^{1/2} t^{\alpha-1/2}$ and noticing that $\alpha \geq \frac{1}{2}$, we see that the set $\{t : Y_t^\alpha(\omega) < 0\}$ is dense near zero for almost all $\omega$. The stopping times $\tau_\delta = \inf\{t > \delta | Y_t^\alpha < -\delta\}$ tend to zero a.s. as $\delta$ goes to zero. Let $\delta$ be small enough such that $A =: \{\tau_\delta < \tau_1\} \in \mathcal{F}_{\tau_\delta}$ has positive measure. Define

$$\tau_\delta^A(\omega) = \begin{cases} \gamma, & \text{when } \omega \notin A, \\ \tau_\delta, & \text{when } \omega \in A, \end{cases}$$

and

$$\tau_1^A(\omega) = \begin{cases} \gamma, & \text{when } \omega \notin A, \\ \tau_1, & \text{when } \omega \in A, \end{cases}$$

where $\gamma$ is a number bigger than the bound of $\tau_1$. Then the process $1_A[Y_t^\alpha - (Y_t^\alpha)^{\tau_\delta^A}]$ admits an equivalent martingale measure $Q$ and $E_Q 1_A Y_{\tau_1^A}^\alpha = E_Q 1_A Y_{\tau_\delta^A}^\alpha$. Since $Y_{\tau_\delta^A}^\alpha < 0$, we have $E_Q 1_A Y_{\tau_1^A}^\alpha < 0$. This shows that $P(Y_{\tau_1}^\alpha < 0) > 0$. Apply the same method and use the upper limit of the law of the iterated logarithm to get $P(Y_{\tau_1}^\alpha > 0) > 0$. We conclude that $Y_t^\alpha$ satisfies the no arbitrage property. □

PROOF OF THEOREM 7. Let $Q$ be the local martingale measure for $S$. Let $\nu_s = \inf\{t > 0 | [S,S]_t > s\}$, then $S_{\nu_s}$ is a Brownian motion under $Q$. We denote it by $B_s$. Then we have $Z_t^\alpha = B_{[S,S]_t} + [S,S]^\alpha$, so $Z_t^\alpha$ is the time changed process of $B_t + t^\alpha$. The result now follows by applying Lemma 6 and Theorem 6. □

**5. Hedging issues.** In this section, we assume the price process $S$ is a continuous semimartingale satisfying condition $(*)$ of Theorem 1. We wish to discuss hedging possibilities.

It is immediately apparent that the restriction of hedging strategies to the class $CC$ greatly reduces the quantity of redundant claims, and essentially all interesting models will be incomplete. Indeed, even in the Brownian paradigm, one cannot have classical completeness without allowing all predictably measurable strategies which are in $L^2(dt \times dP)$, where $dt \times dP$ is understood to be on $[0, T] \times \Omega$. This includes such unrealistic strategies as buying at rational times and selling at irrational times.

Nevertheless, if we are in the Brownian paradigm, we can hope for an *approximate completeness*, in the sense that we can get arbitrarily close to a replicating hedging strategy in an appropriate norm. This idea was developed in a different context in [3], for example (alternatively, see [12]). However, we want to go beyond the usual Itô process framework to include price processes that *normally have arbitrage opportunities*, but do not within our framework of a restricted class of hedging strategies.

Here, *we do not try for maximum generality*, but consider only those claims, which derive from the underlying in a very explicit way, that is, we consider only those contingent claims, which are twice Fréchet continuous functionals of the stock price at time $T$. We further assume the spot interest rate is $r_t \equiv 0$, so we need not worry about the time value of money. Clearly we are not trying for maximum generality here (e.g., in the strict Brownian paradigm, we could replace Fréchet differentiable with Malliavin differentiable), but we are trying only to illustrate what can be done.

Since $S$ is a continuous semimartingale, we can employ the Itô representation formula of Ahn [1], which works for a limited and somewhat special class of contingent claims. We recall Ahn's theorem here for convenience.

THEOREM 8. *Let $f:\mathcal{C}[0,T] \to \mathbb{R}$ be a twice continuously Fréchet differentiable functional at each $x \in \mathcal{C}[0,T]$ with respect to the sup norm. Then for a continuous semimartingale $S$ and $t \in [0,T]$, we have*

$$(5) \quad f(S^t) = f(S^0) + \int_0^t \langle \eta_s, \nabla f(S^s) \rangle \, dS_s + \frac{1}{2} \int_0^t \langle \eta_s \otimes \eta_s, \nabla^2 f(S^s) \rangle \, d[S,S]_s,$$

*where $\eta_s = 1_{[s,T]}$ is an element of the bidual (in the Banach sense) of $\mathcal{C}[0,T]$, and the bracket $\langle \cdot, \cdot \rangle$ is used for dual pairs. Finally, the notation $S^t$ refers to the stopped processes: $S^t_s = S_{s \wedge t}$.*

In Theorem 8 above, we can assume without loss of generality that we are taking the predictable version of the integrands in equation (5).

We need the following lemma.

LEMMA 7. *The space $(CC)$ is dense in $b\mathbb{L}$ in the ucp topology.*

PROOF. By standard results (e.g., Theorem 10 of Chapter II of [16]), we know that simple predictable processes are dense in $b\mathbb{L}$ in ucp, where ucp denotes uniform convergence in probability on compact time sets. Let $(H^n)_{n \geq 1}$ be a sequence of simple predictable processes converging in ucp to $H \in b\mathbb{L}$. Let us fix $\delta > 0$, and choose $\varepsilon > 0$ such that for some $N$ we have if $n \geq N$, then

$$P\left( \sup_{0 \leq t \leq T} |H^n_t - H_t| > \delta \right) \leq \frac{\varepsilon}{2}.$$

Let us now choose and fix an $n_0 \geq N$, and by an abuse of notation, we hereafter refer to $n_0$ simply as $n$. We also suppress the $n\nu$ superscript notation on the jump times. Let $\tau_1 \leq \tau_2 \leq \cdots \leq \tau_n$ be the sequential jump times of $H^n$, and let $\delta(\omega) = \inf_{1 \leq i \leq n}\{\tau_{i+1}(\omega) - \tau_i(\omega)\}$. Since $\tau_{i+1} - \tau_i > 0$ for all $i$ and there are only a finite number of them, there exists a $\delta_0$ such that

$$P\left( \inf_i |\tau^{i+1} - \tau^i| < \delta_0 \right) < \frac{\varepsilon}{2}.$$

Let

$$\Lambda = \left\{ \omega | \inf_i |\tau^{i+1} - \tau^i| < \delta_0 \right\}.$$

Note that we can assume without loss that $P(\Lambda) < \frac{\varepsilon}{2}$ by the strict positivity of the $n < \infty$ random variables $\tau_{i+1} - \tau_i$. Next, let $\nu(\omega) \equiv \inf\{i : |\tau_{i+1} - \tau_i| < \delta_0\}$, and let $\tilde{H}^n_t = H^n_{t \wedge \nu}$.

We need to show that $\tilde{H}^n_t$ is in the space $(CC)$. Since $\tilde{H}^n_t$ changes values only by jumps, we need to show that the jump times of $\tilde{H}^n_t$ are at least $\delta$ apart, for some $\delta > 0$, with $\delta$ nonrandom. Let $\eta_1 < \eta_2 < \cdots < \eta_k$ be the jump times of $\tilde{H}_t^n$. Since $\eta_i = \tau_{i \wedge \nu}$, for $i \leq \nu(\omega)$, the $\eta_i$ exhaust the jumps of $\tilde{H}_t^n$, and since $i \leq \nu(\omega)$, we have that $\eta_i - \eta_{i-1} \geq \delta_0$. We need to show that the random times $\eta_i$ are, in fact, stopping times. To do this, we note that

$$\{i \leq \nu\} = \bigcap_{j=1}^{i} \{\tau_j - \tau_{j-1} \geq \delta_0\} \in \mathcal{F}_{\tau_i}, \tag{6}$$

which implies that

$$\{\eta_i \leq t\} = \{\tau_{i \wedge \nu} \leq t\} = \{\{\tau_i \leq t\} \cap \{i \leq \nu\}\} \bigcup_{j<i} \{\{\tau_j \leq t\} \cap \{\nu = j\}\} \in \mathcal{F}_t.$$
(7)

We conclude that $\tilde{H}_t^n$ is in the class $(CC)$.

To show that it approximates $H^n$, we calculate

$$P\left(\sup_{0 \leq t \leq T} |\tilde{H}^n - H_t| > \delta\right) \leq P\left(\sup_{0 \leq t \leq T} (|\tilde{H}^n - H_{t \wedge \nu}| + |H_{t \wedge \nu} - H_t|) > \delta\right)$$

$$\leq P\left(\sup_{0 \leq t \leq T} |\tilde{H}^n - H_{t \wedge \nu}| > \delta\right)$$

$$+ P\left(\sup_{0 \leq t \leq T} |H_{t \wedge \nu} - H_t| > \delta\right)$$

$$\leq P\left(\sup_{0 \leq t \leq T} |\tilde{H}^n - H_{t \wedge \nu}| > \delta\right) + P(\Lambda)$$

$$\leq \frac{\varepsilon}{2} + \frac{\varepsilon}{2} = \varepsilon.$$

This completes the proof. □

We conclude with our desired result, where we use the semimartingale $\mathcal{H}^2$ norm, as defined in Chapter IV of [16]. This theorem shows that if we first establish a tolerable level of error $\varepsilon$ in the semimartingale $\mathcal{H}^2$ norm, for at least a certain class of contingent claims, we can approximately hedge in class $(CC)$ to within a prescribed error $\varepsilon$.

THEOREM 9. *Let $f : \mathcal{C}[0,T] \to \mathbb{R}$ be a twice continuously Fréchet differentiable functional at each $x \in \mathcal{C}[0,T]$ with respect to the sup norm, and let $C = f(S)$ be a contingent claim whose semimartingale representation given in Theorem 8 is in $\mathcal{H}^2$. Given $\varepsilon > 0$, we can find a hedging strategy $H^n$ in the class $CC$ such that the $\mathcal{H}^2$ norm of $f(S)$, taken on $[0,T]$, is within $\varepsilon$ of the $\mathcal{H}^2$ norm of $\int_0^T H_s^n \, dS_s$.*

PROOF. We know by Theorem 2 of Chapter IV of [16] that the space $b\mathbb{L}$ (bounded, adapted processes with paths which are left continuous with right limits) is dense in bounded predictable processes. By Theorem 4 of Chapter IV of [16], we know that $b\mathbb{L}$ is dense in $b\mathcal{P}$ (the bounded, predictable processes) in the semimartingale $\mathcal{H}^2$ norm, following the definitions given in Chapter IV of [16]. In addition, we recall that simple predictable processes are dense in $b\mathbb{L}$ in the ucp topology, and hence also in the $\mathcal{H}^2$ norm. The result then follows by Lemma 7. □

R. A. JARROW
JOHNONSON GRADUATE SCHOOL OF MANAGEMENT
CORNELL UNIVERSITY
ITHACA, NEW YORK 14853
USA

P. PROTTER
ORIE
CORNELL UNIVERSITY
ITHACA, NEW YORK 14853–381
USA

H. SAYIT
DEPARTMENT OF MATHEMATICS
WORCESTER POLYTECHNIC INSTITUTE
WORCESTER, MASSACHUSETTS 01609
USA
E-MAIL: hs@wpi.edu